\documentclass[12pt]{article}
\usepackage{amssymb,amsfonts,amsmath, psfrag,eepic,colordvi,graphicx,epsfig,ytableau}
\usepackage{amssymb,latexsym,graphics}
\usepackage{MnSymbol}
\usepackage[enableskew]{youngtab}
\usepackage{float}
\usepackage{tikz}
\topskip  
\newcommand{\rmnum}[1]{\romannumeral #1}

 \numberwithin{equation}{section}
\newtheorem{theorem}{Theorem}[section]

\newtheorem{conjecture}[theorem]{Conjecture}

\newtheorem{lemma}[theorem]{Lemma}

\begin{document}
\parskip 6pt

\pagenumbering{arabic}
\def\sof{\hfill\rule{2mm}{2mm}}
\def\ls{\leq}
\def\gs{\geq}
\def\SS{\mathcal S}
\def\qq{{\bold q}}
\def\MM{\mathcal M}
\def\TT{\mathcal T}
\def\EE{\mathcal E}
\def\lsp{\mbox{lsp}}
\def\rsp{\mbox{rsp}}
\def\pf{\noindent {\it Proof.} }
\def\mp{\mbox{pyramid}}
\def\mb{\mbox{block}}
\def\mc{\mbox{cross}}
\def\qed{\hfill \rule{4pt}{7pt}}
\def\block{\hfill \rule{5pt}{5pt}}

\begin{center}
	{\Large\bf   On a conjecture   concerning  the  shuffle-compatible permutation statistics }
\end{center}

\begin{center}
{\small Lihong Yang, Sherry H.F. Yan\footnote{Corresponding author.}  }

 Department of Mathematics\\
Zhejiang Normal University\\
 Jinhua 321004, P.R. China

 hfy@zjnu.cn

\end{center}

\noindent {\bf Abstract.} The  notion of shuffle-compatible permutation statistics  was implicit in Stanley's work on P-partitions and was first
explicitly studied by Gessel and Zhuang.  The aim of this paper is to prove that the triple  ${\rm (udr, pk, des)}$ is shuffle-compatible as conjectured by Gessel and Zhuang, where ${\rm udr}$  denotes the number of up-down runs,   ${\rm pk}$ denotes the peak number,   and ${\rm des}$ denotes the descent number.  This is accomplished by establishing an ${\rm (udr, pk, des)}$-preserving bijection in the spirit of  Baker-Jarvis and  Sagan's   bijective proofs of   shuffle-compatibility property of   permutation statistics. As an application, our bijection also enables us to prove that
the pair $({\rm cpk}, {\rm cdes})$ is cyclic  shuffle-compatible, where ${\rm cpk}$ denotes the cyclic peak number   and ${\rm cdes}$ denotes the cyclic descent number.



%===========================================================================

\section{Introduction}
Let $\mathbb{P}$ denote the set of  all  positive
integers. To denote the
cardinality of a set  $U$, we use $|U|$.  For  $U\subset \mathbb{P}$ with $|U|=n$, a permutation of $U$ is a linear order $\pi=\pi_1\pi_2\ldots \pi_n$ of the elements of $U$.
   Denote by $L(U)$ the set of all permutations of $U$. The {\em length} of a permutation $\pi$ is the cardinality of its underlying set, i.e. $|U|$, which  is denoted
  by $|\pi|$. Permutations  have  been  extensively  studied  over the last  decades.   For a thorough summary of the current status of research, see B\'{o}na's   book \cite{Bona}.

 The three classical examples of permutation statistics are the descent set ${\rm Des}$, the descent
number ${\rm des}$, and the major index ${\rm maj}$. For $\pi\in L(U)$ with $|U|=n$,   define
$$
 {\rm Des}(\pi) = \{i: \pi_i>\pi_{i+1}, 1\leq i\leq n-1\},
 $$
 $$
{\rm des}(\pi) = |{\rm Des}(\pi)|,
$$
and
$$
{\rm maj}(\pi) = \sum_{i\in {\rm Des}(\pi)}i.
 $$

 A statistic  ${\rm st}$ is said to be a {\em descent statistic} if  ${\rm Des}(\pi)={\rm Des}(\sigma)$ implies that ${\rm st}(\pi)={\rm st}(\sigma)$ for any two permutations $\pi$ and $\sigma$. Clearly, the statistics ${\rm Des}$, ${\rm des}$ and ${\rm maj}$ are descent statistics.   For $\pi\in L(U)$ with $|U|=n$,   the {\em peak set} of $\pi$, denoted by ${\rm Pk}(\pi)$, is defined to be
 $$
 {\rm Pk}(\pi)=\{i: \pi_{i-1}<\pi_i>\pi_{i+1}, 2\leq i\leq n-1\}.
 $$
 The {\em peak number} of  $\pi$, denoted by ${\rm pk}(\pi)$, is defined to be the cardinality of ${\rm Pk}(\pi)$.   The {\em  exterior peak number} of  $\pi$, denoted by ${\rm epk}(\pi)$, is defined to be the peak number of the permutation $0\pi0$.
 A {\em monotone factor} of a permutation  is a factor that is either strictly increasing or
strictly decreasing. A {\em birun} is a maximal monotone factor. An {\em updown run} is a birun
of $0\pi$. The number of biruns and  updown runs of $\pi$ are  denoted ${\rm bir}(\pi)$ and  ${\rm udr}(\pi)$,  respectively.

For any two permutations $\pi\in L(U)$ and $\sigma\in L(V)$ with $U\bigcap V=\emptyset$, we say
 that the permutation $\tau\in L(U\bigcup V)$ is a {\em shuffle} of $\pi$ and $\sigma$  if both $\pi$ and $\sigma$ are subsequences of $\tau$.  Denote by $S(\pi, \sigma)$ the set of shuffles of $\pi$ and $\sigma$. For example, $S(31, 24)=\{3124, 3241, 2431, 3214, 2341, 2314\}$. A permutation  statistic ${\rm st}$ is said to be {\em  shuffle-compatible} if for any   permutations $\pi$ and $\sigma$ with disjoint underlying  sets, the multiset
$\{ {\rm st}(\tau) : \tau \in  S(\pi, \sigma) \}$,  which   encodes the distribution of the statistic ${\rm st}$ over
shuffles of $\pi$ and $
\sigma$, depends only on ${\rm st}(\pi)$, ${\rm st}(\sigma)$, $|\pi|$  and $|\sigma|$. For our convenience, we simply write $st(S(\pi, \sigma))$  for the multiset $\{ {\rm st}(\tau) : \tau \in  S(\pi, \sigma) \}$.  For instance, ${\rm des}(S(31, 24))=\{1^3,2^3\}$.
We say that the permutation  statistic ${\rm st}$ has {shuffle-compatibility  property } if ${\rm st}$ is shuffle-compatible.

For nonnegative integer $n$, let
$$
[n]_q=1+q+q^2+\ldots +q^{n-1}
$$ and
$$
[n]_q!=[1]_q[2]_q\ldots [n]_q
$$
 where $q$   is a variable. For $0\leq k\leq n$, let
 $$
 {n\choose k}_q={[n]_q!\over [k]_q![n-k]_q!}.
 $$
By utilizing P-partitions, Stanley \cite{Stanley} proved that for any two permutations $
\pi$ and $\sigma$ with disjoint  underlying sets,
 \begin{equation}\label{seq2}
\begin{array}{lll}
\sum_{\tau\in S_k(\pi, \sigma)}q^{{\rm maj}(\tau)}&=&q^{{\rm maj}(\pi)+{\rm maj}(\sigma)+(k-{\rm des}(\pi))(k-{\rm des}(\sigma))}{|\pi|-{\rm des}(\pi)+{\rm des}(\sigma)\choose k-{\rm des}(\pi) }_{q}\\
& & \times {|\sigma|-{\rm des}(\sigma)+{\rm des}(\pi)\choose k-{\rm des}(\sigma) }_{q}
\end{array}
\end{equation}
where $S_k(\pi, \sigma)=\{\tau: \tau\in S(\pi, \sigma), {\rm des}(\tau)=k\}$.
 The bijective proofs of (\ref{seq2})  have been given by Goulden \cite{Goulden},
Stadler \cite{Stadler}, Ji and Zhang \cite{Ji},  respectively.    Novick \cite{Novick} provided a bijective proof  of  the following formula  due to Garsia and Gessel \cite{Gar}:
\begin{equation}\label{seq1}
\sum_{\tau\in S(\pi, \sigma)}q^{{\rm maj}(\tau)}=q^{{\rm maj}(\pi)+{\rm maj}(\sigma)}{|\pi|+|\sigma|\choose |\pi|}_{q}
\end{equation}
where
$
\pi$ and $\sigma$ are permutations with disjoint underlying sets. Very recently,  Ji and Zhang  \cite{Ji2} derived  a cyclic analogue of (\ref{seq2}). 
Formulae (\ref{seq2}) and (\ref{seq1}) imply  that the statistics ${\rm maj}$ and $({\rm maj}, {\rm des})$ are shuffle-compatible.

By using noncommutative symmetric functions,
quasisymmetric functions, and variants of quasisymmetric functions,
Gessel and Zhuang  \cite{Gessel}   further  investigated  the shuffle-compatibility property of  permutation statistics and proved that many permutation statistics do have this property.  They also posed several conjectures concerning the shuffle-compatibility of permutation statistics.   Some of these conjectures were then confirmed  by Grinberg \cite{Grinberg} and O$\breve{g}$uz \cite{Oguz}.
Recently,   Baker-Jarvis and  Sagan \cite{Baker} presented a bijective approach to  deal with the
  shuffle compatibility of permutations statistics. As an application, Baker-Jarvis and  Sagan \cite{Baker} proved that the pair $({\rm udr}, {\rm pk}) $ is shuffle-compatible as conjectured  by Gessel and Zhuang  \cite{Gessel}. 

The main objective of this paper is to prove  the following conjecture posed by Gessel and Zhuang  \cite{Gessel}.
\begin{conjecture}{\upshape   (See \cite{Gessel} , conjecture 6.7 )}\label{con1}
The   triple  $ (\rm udr,\rm pk,\rm des)  $	is  shuffle-compatible.
\end{conjecture}

In \cite{Baker}, Baker-Jarvis and  Sagan remarked that their bijection for proving the shuffle compatibility of the  statistic $( \rm udr,\rm pk)$ does not preserve the statistic ${\rm des}$ and posed an  open problem of finding a bijective proof of the the shuffle compatibility of the  statistic $( \rm udr,\rm pk, \rm des)$ (see \cite{Baker}, Question 7.1). In this paper, we  aim to provide such a bijective proof in the spirit of  Baker-Jarvis and  Sagan's   bijective proofs of shuffle-compatibility property of   permutation statistics.

Recently, Adin, Gessel, Reiner and Roichman \cite{Adin} introduced a cyclic version of quasisymmetric functions
with a corresponding cyclic shuffle operation. A cyclic permutation $[\pi]$ of $U$ can be viewed
as an equivalence class of linear permutations  $\pi=\pi_1\pi_2\ldots \pi_n$ of  $U$ under the cyclic equivalence relation
$\pi_1\pi_2\ldots \pi_n \sim \pi_i\ldots \pi_n \pi_1\pi_2\ldots \pi_{i-1}$ for all $2\leq i\leq n$.  For example
$$
[1243]=\{1243, 2431, 4312, 3124\}
$$
is a cyclic permutation of $U=[4]$. Denote by $C(U)$ the set of all cyclic permutations of $U$.  Let $\pi_{\ell}$ be the smallest element of $U$, then the linear permutation $\pi_{\ell}\pi_{\ell+1}\ldots \pi_n\pi_1\pi_2\ldots \pi_{\ell-1}
$ is called the {\em representative } of the cyclic permutation $[\pi]$. For the example above, $1243$ is the representative of the cyclic permutation $[1243]$. Here and in the sequel, we use the representative
to represent each cyclic permutation.
   For example, for $U=[4]$,   the elements of $C(U)$ are listed as follows:
$$
 [1234], [1243], [1324], [1342], [1423], [1432].
$$
For a linear permutation $\pi=\pi_1\pi_2\ldots \pi_n$, define  the {\em cyclic descent set  } and the {\em cyclic descent number } of $\pi$   to be
$$
{\rm cDes}(\pi)=\{i\mid \pi_i>\pi_{i+1}\}
$$
and
$$
{\rm cdes}(\pi)=|{\rm cDes}(\pi)|
$$
with the convention $\pi_{n+1}=\pi_1$. Similarly, the  {\em cyclic peak set  }  and the {\em cyclic peak number }of $\pi$ are  defined to be
$$
{\rm cPk}(\pi)=\{i\mid \pi_{i-1}<\pi_i>\pi_{i+1}\}
$$
and
$$
{\rm cpk}(\pi)=|{\rm cPk}(\pi)|
$$
with the convention $\pi_{n+1}=\pi_1$ and $\pi_0=\pi_{n}$.
For example, let $\pi=4218596$.  We have
		$
		{\rm cDes}(\pi)=\{1,2,4,6,7\}, {\rm cdes}(\pi)=5, {\rm cPk}(\pi)=\{4,6\}, {\rm cpk}(\pi)=2.
		$

For a  cyclic permutation $[\pi]$, 		define  the {\em cyclic descent set  } and  {\em cyclic peak set} of $\pi$   to be
$$
{\rm cDes}([\pi])=\{\{{\rm cDes}(\sigma)\}\mid \sigma\in [\pi]\},
$$
and
$$
{\rm cPk}([\pi])=\{\{{\rm cPk}(\sigma)\}\mid \sigma\in [\pi]\}.
$$
Define the   {\em cyclic descent number } and  {\em cyclic peak number } of $\pi$ to be
$$
{\rm cdes}([\pi])={\rm cdes}(\pi)
$$
and
$$
{\rm cpk}([\pi])={\rm cpk}(\pi).
$$

 For any two cyclic  permutations $[\pi]\in C(U)$ and $[\sigma]\in C(V)$ with $U\bigcap V=\emptyset$, we say
 that the cyclic permutation $[\tau]\in C(U\bigcup V)$ is a {\em cyclic shuffle} of $[\pi]$ and $[\sigma]$  if both $[\pi]$ and $[\sigma]$ are circular subsequences of $[\tau]$.  Denote by $cS([\pi], [\sigma])$ the set of cyclic  shuffles of $[\pi]$ and $[\sigma]$. For example, let $[\pi]=[13 ]$ and $[\sigma]=[24]$. We have
		$$
 cS([\pi], [\sigma])=\{[1423], [1342], [1432], [1234], [1324], [1243]\}.
		$$
For a cyclic permutation statistic ${\rm cst}$, define  ${\rm cst}(cS([\pi], [\sigma]))$ to be the multiset $\{ {\rm cst}([\tau]) : [\tau ]\in  cS([\pi], [\sigma]) \}$. Continuing with the above example, we have
	$$
 {\rm cdes}(cS([\pi], [\sigma]))=\{1,2^4, 3\}
		$$  and $${\rm cpk} (cS([\pi], [\sigma]))=\{1^4, 2^2\}.$$
A cyclic permutation  statistic ${\rm cst}$ is said to be {\em  cyclic shuffle-compatible} if for any  cyclic  permutations $[\pi]$ and $[\sigma]$ with disjoint underlying  sets, the multiset ${\rm cst}(cS([\pi], [\sigma]))$  depends only on ${\rm cst}([\pi])$, ${\rm cst}([\sigma])$, $| \pi |$  and $| \sigma |$.

Very recently, Domagalski, Liang, Minnich, Sagan, Schmidt and Sietsema \cite{Doma} derived  the following cyclic
  shuffle compatibility results.
\begin{theorem}{\upshape   (See \cite{Doma} ,  Theorem  1.2 )}
The statistics
$$
{\rm cDes}, {\rm cPk }, {\rm cdes}, {\rm cpk}
$$
are cyclic  shuffle-compatible.
\end{theorem}

Gessel and Zhuang \cite{Gessel} proved that the pair $({\rm des}, {\rm pk})$ is    shuffle-compatible.  In this paper, we will prove the following cyclic analogue of Gessel and Zhuang's result.
\begin{theorem}\label{cthm}
The pair $({\rm cpk}, {\rm cdes})$ is cyclic  shuffle-compatible.
\end{theorem}

\section{Proof of Conjecture \ref{con1}}
This section is devoted to the bijective proof of Conjecture \ref{con1}. To this end, we need to recall the  following  two lemmas due to  Baker-Jarvis and  Sagan \cite{Baker}.

\begin{lemma}{\upshape   (See \cite{Baker} ,  Theorem 4.2 )}\label{lem1}
The statistic ${\rm Des}$ is shuffle-compatible.
\end{lemma}

 For $m,n\geq 1$, let $[n]=\{1,2,\ldots, n\}$ and $[n]+m=\{n+i: 1\leq i\leq m\}$.

\begin{lemma}{\upshape   (See \cite{Baker} ,  Corollary 3.2 )}\label{lem2}
Suppose that  ${\rm st}$ is a descent statistic. The following are equivalent.
\begin{itemize}
\item[(a)]  The statistic ${\rm st}$ is shuffle-compatible.
\item[(b)] If  ${\rm st} (\pi)={\rm st}(\pi')$ where  $\pi, \pi'\in L([n])$, and $\sigma\in L([n]+m)$ for some $m, n\geq 1$, then $
{\rm st}(S(\pi, \sigma) )={\rm st}(S(\pi', \sigma)).
$

\end{itemize}

\end{lemma}

 For a   permutation $\pi\in L(U)$ with $k$ biruns, the {\em type} of $\pi$, denoted by ${\rm type}(\pi)$, is defined to be $ (t_1, t_2, \ldots, t_k)$, where $t_i$ denotes the length of the $i$-th  birun  (counting from left to right). For example,  ${\rm type}(6534792)=(3, 4,2)$. For a permutation $\pi=\pi_1\pi_2\ldots \pi_n$, define $\chi^+(\pi)$ to be $1$ if $\pi_1>\pi_2$ and  to be $0$ otherwise. Similarly,  we define $\chi^-(\pi)$ to be $1$ if $\pi_{n-1}<\pi_n$ and  to be $0$ otherwise.  One can easily check that
 \begin{equation}\label{pkudr}
 {\rm udr}(\pi)=\left\{
 \begin{array}{ll}
 2{\rm pk}(\pi)& \mbox {if}\,\,  \chi^{+}(\pi)=\chi^{-}(\pi)=0, \\
   2{\rm pk}(\pi)+1& \mbox{if}\,\,  \chi^{+}(\pi)=0, \chi^{-}(\pi)=1, \\
    2{\rm pk}(\pi)+2& \mbox{if}\,\,  \chi^{+}(\pi)=1, \chi^{-}(\pi)=0, \\
   2{\rm pk}(\pi)+3& \mbox{if}\,\,  \chi^{+}(\pi)= \chi^{-}(\pi)=1. \\
 \end{array}
 \right.
 \end{equation}

Let $\pi\in L([n])$ be a permutation  with  ${\rm type}(\pi)=(t_1, t_2, \ldots, t_k)$ such that $t_{\ell}\geq 3$ for some $\ell\geq 3$.
Define $\Omega_{\ell}(\pi)$ to be the set of  permutations $\pi'\in L([n])$ with $\chi^{+}(\pi')=\chi^{+}(\pi)$  and ${\rm type}(\pi')=(t'_1, t'_2,  \ldots t'_k)$  where
$$ t'_i=\left\{
\begin{array}{ll}
t_{i}+1 & \mbox{if} \,\, i=\ell-2,\\
t_{i}-1 & \mbox{if} \,\, i=\ell,\\
t_{i}& \mbox{otherwise.}
\end{array}
\right.
$$
One can easily check that for any $
\pi'\in  \Omega_{\ell}(\pi)$, we have ${\rm (udr, pk, des)}\pi={\rm (udr, pk, des)}\pi'$.

In order to prove Conjecture \ref{con1}, we define  four disjoint  canonical sets as follows.
Define
$$
\Pi^{(1)}_{n,k,d}=\{\pi\in L([n]): \chi^{+}(\pi)=0,  {\rm type}(\pi)=(t_1, t_2, \ldots, t_{2k})\}
$$
where $t_1=n-d-k+1$, $t_2=d-k+2$,  and $t_i=2$ for $2<i\leq 2k$. For example, we have
 $
 \pi=25796431(10)8 \in \Pi^{(1)}_{10,2,5}$ with   $ {\rm type}(\pi)=(4,5,2,2)
 $

Define
$$
\Pi^{(2)}_{n,k,d}=\{\pi\in L([n]): \chi^{+}(\pi)=0,  {\rm type}(\pi)=(t_1, t_2, \ldots, t_{2k+1})\}
$$
where $t_1=n-d-k$, $t_2=d-k+2$,  and $t_i=2$ for $2<i\leq 2k+1$. For example, we have
$
\pi=2796431(10)58 \in \Pi^{(2)}_{10,2,5} $ with $ {\rm type}(\pi)=(3,5,2,2,2)
$

Define
$$
\Pi^{(3)}_{n,k,d}=\{\pi\in L([n]): \chi^{+}(\pi)=1,  {\rm type}(\pi)=(t_1, t_2, \ldots, t_{2k+1})\}
$$
where $t_1=d-k+1$, $t_2=n-d-k+1$,  and $t_i=2$ for $2<i\leq 2k+1$. For example, we have
$
\pi=964123(10)785 \in \Pi^{(3)}_{10,2,5}$ with ${\rm type}(\pi)=(4,4,2,2,2)
$

Define
$$
\Pi^{(4)}_{n,k,d}=\{\pi\in L([n]): \chi^{+}(\pi)=1,  {\rm type}(\pi)=(t_1, t_2, \ldots, t_{2k+2})\}
$$
where $t_1=d-k+1$, $t_2=n-d-k$,  and $t_i=2$ for $2<i\leq 2k+2$. For example, we have
$
\pi=96412(10)3857 \in \Pi^{(4)}_{10,2,5}$ with $ {\rm type}(\pi)=(4,3,2,2,2,2).
$
Let
 $$
 \Pi_{n,k,d}= \Pi^{(1)}_{n,k,d}\bigcup \Pi^{(2)}_{n,k,d}\bigcup \Pi^{(3)}_{n,k,d}\bigcup \Pi^{(4)}_{n,k,d}.
$$

By (\ref{pkudr}),  one  can deduce the following result.
\begin{lemma}\label{lempk}
 For any permutation $
\pi\in \Pi_{n, k,d}$, we have
$$
({\rm udr, pk})\pi=\left\{
\begin{array}{ll}
(2k, k)& \mbox{if }\,\, \pi\in \Pi^{(1)}_{n,k,d},\\
 (2k+1, k)& \mbox{if }\,\, \pi\in \Pi^{(2)}_{n,k,d},\\
 (2k+2, k)& \mbox{if }\,\, \pi\in \Pi^{(3)}_{n,k,d},\\
 (2k+3, k)& \mbox{if }\,\, \pi\in \Pi^{(4)}_{n,k,d}.\\
\end{array}
\right.
$$
\end{lemma}

The following theorem  will play an essential role in the proof of Conjecture \ref{con1}.

\begin{theorem}\label{th1}
Let $\pi\in L([n])$ be a permutation  with $({\rm pk, des }) \pi =(k,d)$ and let $\sigma\in L([n]+m)$ for some $m,n\geq 1$ and $k,d\geq 0$.  The following statements hold.
\begin{itemize}
\item[\upshape (\rmnum{1}) ] If ${\rm type}(\pi)=(t_1, t_2, \ldots, t_{2k})$, $\chi^{+}(\pi)=0$,  and  $\pi\notin  \Pi^{(1)}_{n,k,d}$, then there exists  a permutation $\pi'\in  \Pi^{(1)}_{n,k,d}$ such that $$
 {(\rm udr, pk, des)}\pi= {(\rm udr, pk, des)}\pi'
 $$ and $$ {(\rm udr, pk, des)}S(\pi, \sigma)= {(\rm udr, pk, des)}S(\pi', \sigma).  $$
 \item[\upshape (\rmnum{2}) ] If ${\rm type}(\pi)=(t_1, t_2, \ldots, t_{2k+1})$, $\chi^{+}(\pi)=0$,  and  $\pi\notin  \Pi^{(2)}_{n,k,d}$, then there exists  a permutation $\pi'\in  \Pi^{(2)}_{n,k,d}$ such that $$
 {(\rm udr, pk, des)}\pi= {(\rm udr, pk, des)}\pi'
 $$ and $$ {(\rm udr, pk, des)}S(\pi, \sigma)= {(\rm udr, pk, des)}S(\pi', \sigma).  $$.
 \item[\upshape (\rmnum{3}) ] If ${\rm type}(\pi)=(t_1, t_2, \ldots, t_{2k+1})$, $\chi^{+}(\pi)=1$,  and  $\pi\notin  \Pi^{(3)}_{n,k,d}$, then there exists  a permutation $\pi'\in  \Pi^{(3)}_{n,k,d}$ such that $$
 {(\rm udr, pk, des)}\pi= {(\rm udr, pk, des)}\pi'
 $$ and $$ {(\rm udr, pk, des)}S(\pi, \sigma)= {(\rm udr, pk, des)}S(\pi', \sigma).  $$
 \item[\upshape (\rmnum{4}) ] If ${\rm type}(\pi)=(t_1, t_2, \ldots, t_{2k+2})$, $\chi^{+}(\pi)=1$,  and  $\pi\notin  \Pi^{(4)}_{n,k,d}$, then there exists  a permutation $\pi'\in  \Pi^{(4)}_{n,k,d}$ such that $$
 {(\rm udr, pk, des)}\pi= {(\rm udr, pk, des)}\pi'
 $$ and $$ {(\rm udr, pk, des)}S(\pi, \sigma)= {(\rm udr, pk, des)}S(\pi', \sigma).  $$

\end{itemize}
 \end{theorem}

Before we prove Theorem  \ref{th1} , we need the following lemma.
\begin{lemma}\label{mlem1}
Let $\pi\in L([n])$ be a permutation  with  ${\rm type}(\pi)=(t_1, t_2, \ldots, t_k)$ such that $t_{\ell}\geq 3$ for some $\ell\geq 3$ and let $\sigma\in L([n]+m)$ for some $m,n\geq 1$. Then  there exists an ${\rm (udr, pk, des)}$-preserving  bijection $\phi_{\ell}: S(\pi, \sigma)\longrightarrow S(\pi', \sigma)$ for any permutation $\pi'\in \Omega_{\ell}(\pi)$.

\end{lemma}
\pf    Let $\tau=\tau_1\tau_2\ldots \tau_{n+m}\in S(\pi, \sigma)$. If the $\ell$-th birun is increasing (resp. decreasing), then let   $\pi_j$ and $\pi_{j+1}$  be the first (resp. last) two entries of $\ell$-th birun of $\pi$ and let $\pi_i$ be the first (resp. last)  entry of the $(\ell-2)$-th birun of $\pi$.  Then $\tau$ can be uniquely factored as $\tau^a \tau^{b}\tau^{c}$, where $\tau^{b}$ is the  subsequence of $\tau$ between  $\pi_i$ and $\pi_{j+1}$ including $\pi_i$ and $\pi_{j+1}$.
Then $\tau^{b}$ can be further decomposed  as
$$
\pi_{i}\sigma^{(1)}\pi_{i+1}\sigma^{(2)}\ldots \pi_j \sigma^{(j-i+1)}\pi_{j+1},
$$
where   $\sigma^{(s)}$ is a (possibly empty) subsequence of $\tau$ and all the entries of $\sigma^{(s)}$ belong  to $\sigma$ for all $1\leq s\leq j-i+1$.   Now we proceed to construct $\phi_{\ell}(\tau)$ by distinguishing the following two cases.

\noindent {\bf  Case 1.}  $\sigma^{(j-i+1)}=\emptyset$.\\
Define $\phi_{\ell}(\tau)$ to be the permutation $\theta^a \theta^{b}\theta^{c}$, where $\theta^{a}$   (resp. $\theta^c$ ) is the permutation obtained from $\tau^a$ (resp. $\tau^c$ ) by replacing each element $\pi_k$ by $\pi'_k$ for $1\leq k<i$ (resp. $j+1<k\leq n$)  and
    $$
    \theta^{b}=\pi'_{i}\pi'_{i+1}\sigma^{(1)}\pi'_{i+2}\sigma^{(2)}\ldots \pi'_{j}\sigma^{(j-i)} \pi'_{j+1}.
    $$

   For example,  let $\ell=4$,   $\pi=6351274\in L([7])$ and $\sigma=(11)89(10) \in L([7]+4)$. Then  $\tau=6 3 {\bf (11) 8}5{\bf 9}127{\bf (10)}4\in S(\pi,\sigma)$ and   $\pi'=6145273\in \Omega_{4}(\pi)$.   Then $\tau$ can be decomposed as $\tau^{a}\tau^{b}\tau^{c}$ as illustrated in Figure  \ref{f1}.
   Clearly,  $\tau^{b}$ can be further  decomposed  as $3\sigma^{(1)}5\sigma^{(2)}1\sigma^{(3)}2$ where $\sigma^{(1)}={\bf (11)8}$, $\sigma^{(2)}={\bf 9}$ and $\sigma^{(3)}=\emptyset$.
   By applying the map $\phi_{4}$ to $\tau$,   we obtain $\phi_{4}(\tau)=\theta^{a}\theta^{b}\theta^{c}$ as shown in Figure \ref{f1}, where $\theta^{a}=6$, $\theta^{b}=14{\bf (11)8}5{\bf 9}2$ and $\theta^{c}=7{\bf (10)}3$.

\begin{figure}[H]
\begin{center}
	\begin{tikzpicture}
     \draw[step=0.4cm,gray,very thin] (-2,-2) grid(2.8,2.4);

     \filldraw[fill=black](-1.6,0.4)circle(0.1);
     \filldraw[fill=black](-1.2,-0.8)circle(0.1);
     \filldraw[fill=red](-0.8,2.4)circle(0.1);
     \filldraw[fill=red](-0.4,1.2)circle(0.1);
     \filldraw[fill=black](0,0)circle(0.1);
     \filldraw[fill=red](0.4,1.6)circle(0.1);
     \filldraw[fill=black](0.8,-1.6)circle(0.1);
     \filldraw[fill=black](1.2,-1.2)circle(0.1);
     \filldraw[fill=black](1.6,0.8)circle(0.1);
     \filldraw[fill=red](2,2)circle(0.1);
     \filldraw[fill=black](2.4,-0.4)circle(0.1);

     \draw[dotted][very thick](-1.6,0.4)--(-1.2,-0.8);
     \draw[dotted][very thick](-0.8,2.4)--(-1.2,-0.8);
     \draw[dotted][very thick](-0.8,2.4)--(-0.4,1.2);
     \draw[dotted][very thick](0,0)--(-0.4,1.2);
     \draw[dotted][very thick](0,0)--(0.4,1.6);
     \draw[dotted][very thick](0.8,-1.6)--(0.4,1.6);
     \draw[dotted][very thick](0.8,-1.6)--(1.2,-1.2);
     \draw[dotted][very thick](1.6,0.8)--(1.2,-1.2);
     \draw[dotted][very thick](1.6,0.8)--(2,2);
     \draw[dotted][very thick](2.4,-0.4)--(2,2);

     \draw[blue][ultra thick](-1.2,-2)--(-1.2,2.4);
     \draw[blue][ultra thick](1.2,-2)--(1.2,2.4);

     \coordinate [label=above:$6$] (x) at (-1.6,-2.8);
     \coordinate [label=above:$3$] (x) at (-1.2,-2.8);
     \coordinate [label=above:$11$] (x) at (-0.8,-2.8);
     \coordinate [label=above:$8$] (x) at (-0.4,-2.8);
     \coordinate [label=above:$5$] (x) at (0,-2.8);
     \coordinate [label=above:$9$] (x) at (0.4,-2.8);
     \coordinate [label=above:$1$] (x) at (0.8,-2.8);
     \coordinate [label=above:$2$] (x) at (1.2,-2.8);
     \coordinate [label=above:$7$] (x) at (1.6,-2.8);
     \coordinate [label=above:$10$] (x) at (2,-2.8);
     \coordinate [label=above:$4$] (x) at (2.4,-2.8);

     \coordinate [label=above:$\tau^{a}$] (x) at (-1.5,2.4);
     \coordinate [label=above:$\tau^{b}$] (x) at (0.1,2.4);
     \coordinate [label=above:$\tau^{c}$] (x) at (2.1,2.4);

     \draw[->] (3.6,0)--(5.2,0);\put(40,1.2){$\phi_{4}(\tau)$};

     \draw[step=0.4cm,gray,very thin] (5.999,-2) grid(10.8,2.4);

     \filldraw[fill=black](6.4,0.4)circle(0.1);
     \filldraw[fill=black](6.8,-1.6)circle(0.1);
     \filldraw[fill=black](7.2,-0.4)circle(0.1);
     \filldraw[fill=red](7.6,2.4)circle(0.1);
     \filldraw[fill=red](8,1.2)circle(0.1);
     \filldraw[fill=black](8.4,0)circle(0.1);
     \filldraw[fill=red](8.8,1.6)circle(0.1);
     \filldraw[fill=black](9.2,-1.2)circle(0.1);
     \filldraw[fill=black](9.6,0.8)circle(0.1);
     \filldraw[fill=red](10,2)circle(0.1);
     \filldraw[fill=black](10.4,-0.8)circle(0.1);

     \draw[dotted][very thick](6.4,0.4)--(6.8,-1.6);
     \draw[dotted][very thick](6.8,-1.6)--(7.2,-0.4);
     \draw[dotted][very thick](7.2,-0.4)--(7.6,2.4);
     \draw[dotted][very thick](7.6,2.4)--(8,1.2);
     \draw[dotted][very thick](8,1.2)--(8.4,0);
     \draw[dotted][very thick](8.4,0)--(8.8,1.6);
     \draw[dotted][very thick](8.8,1.6)--(9.2,-1.2);
     \draw[dotted][very thick](9.2,-1.2)--(9.6,0.8);
     \draw[dotted][very thick](9.6,0.8)--(10,2);
     \draw[dotted][very thick](10,2)--(10.4,-0.8);

     \draw[blue][ultra thick](6.8,-2)--(6.8,2.4);
     \draw[blue][ultra thick](9.2,-2)--(9.2,2.4);

     \coordinate [label=above:$6$] (x) at (6.4,-2.8);
     \coordinate [label=above:$1$] (x) at (6.8,-2.8);
     \coordinate [label=above:$4$] (x) at (7.2,-2.8);
     \coordinate [label=above:$11$] (x) at (7.6,-2.8);
     \coordinate [label=above:$8$] (x) at (8,-2.8);
     \coordinate [label=above:$5$] (x) at (8.4,-2.8);
     \coordinate [label=above:$9$] (x) at (8.8,-2.8);
     \coordinate [label=above:$2$] (x) at (9.2,-2.8);
     \coordinate [label=above:$7$] (x) at (9.6,-2.8);
     \coordinate [label=above:$10$] (x) at (10,-2.8);
     \coordinate [label=above:$3$] (x) at (10.4,-2.8);

     \coordinate [label=above:$\theta^{a}$] (x) at (6.5,2.4);
     \coordinate [label=above:$\theta^{b}$] (x) at (8.2,2.4);
     \coordinate [label=above:$\theta^{c}$] (x) at (10.2,2.4);

	\end{tikzpicture}
\end{center}
\caption{ An example of  Case 1.}\label{f1}
\end{figure}

  \noindent {\bf  Case 2.}  $\sigma^{(j-i+1)}\neq \emptyset$. \\
  Suppose that $\sigma^{(s)}\neq \emptyset$  if and only if  $s\in \{s_1, s_2, \ldots, s_p\}$ with $1\leq s_1<s_2<\ldots< s_p=j-i+1$.
  Define $\phi_{\ell}(\tau)$ to be the permutation $\theta^a \theta^{b}\theta^{c}$, where $\theta^{a}$   (resp. $\theta^c$ ) is the permutation obtained from $\tau^a$ (resp. $\tau^c$ ) by replacing each element $\pi_k$ with $\pi'_k$  for $1\leq k<i$ (resp. $j+1<k\leq n$)  and
    $
    \theta^{b}
    $ is obtained from $\tau^b$ by replacing each $\pi_k$ with  $\pi'_{k+1}$ for $i\leq k\leq j$,   replacing each $\sigma^{(s_q)}$ by  $\sigma^{(s_{q+1})}$ for $1\leq q\leq  p-1$,  and inserting the subsequence  $\pi'_i\sigma^{(s_1)}$ immediately to the left of  $\pi'_{i+1}$.

     For example, let $\ell=3$, $\pi=7426315\in L([7])$ and $\sigma=(11)8(10)9(12) \in L([7]+5)$. Then $\tau={\bf (11)}74{\bf 8}2{\bf (10)}63{\bf 9}{\bf (12)}15\in S(\pi,\sigma)$ and $\pi'=7432615 \in \Omega_{3}(\pi)$. Figure \ref{f4} illustrates the decomposition of    $\tau$, where $\tau^{a}={\bf (11)}74{\bf 8}$, $\tau^{b}=2{\bf (10)}63{\bf 9}{\bf (12)}1$ and $\tau^{c}=5$. Clearly, $\tau^{b}$ can be further decomposed as $ 2\sigma^{(1)}6\sigma^{(2)}3\sigma^{(3)}1$ where $\sigma^{(1)}={\bf (10)}$, $\sigma^{(2)}=\emptyset$,  and $\sigma^{(3)}={\bf 9}{\bf (12)}$.   By applying the map $\phi_{3}$ to $\tau$,   we obtain $\phi_{3}(\tau)=\theta^{a}\theta^{b}\theta^{c}$ as shown in Figure \ref{f4},  where $\theta^{a}={\bf (11)}74{\bf 8}$, $\theta^{b}=3{\bf (10)}2{\bf 9}{\bf (12)}61$ and $\theta^c=5$.

    \begin{figure}[H]
   	\begin{center}
   		\begin{tikzpicture}
   			\draw[step=0.4cm,gray,very thin] (-2,-2) grid(3.2,2.8);
   			
   			\filldraw[fill=red](-1.6,2.4)circle(0.1);
   			\filldraw[fill=black](-1.2,0.8)circle(0.1);
   			\filldraw[fill=black](-0.8,-0.4)circle(0.1);
   			\filldraw[fill=red](-0.4,1.2)circle(0.1);
   			\filldraw[fill=black](0,-1.2)circle(0.1);
   			\filldraw[fill=red](0.4,2)circle(0.1);
   			\filldraw[fill=black](0.8,0.4)circle(0.1);
   			\filldraw[fill=black](1.2,-0.8)circle(0.1);
   			\filldraw[fill=red](1.6,1.6)circle(0.1);
   			\filldraw[fill=red](2,2.8)circle(0.1);
   			\filldraw[fill=black](2.4,-1.6)circle(0.1);
   			\filldraw[fill=black](2.8,0)circle(0.1);
   			
   			\draw[dotted][very thick](-1.6,2.4)--(-1.2,0.8);
   			\draw[dotted][very thick](-1.2,0.8)--(-0.8,-0.4);
   			\draw[dotted][very thick](-0.8,-0.4)--(-0.4,1.2);
   			\draw[dotted][very thick](-0.4,1.2)--(0,-1.2);
   			\draw[dotted][very thick](0,-1.2)--(0.4,2);
   			\draw[dotted][very thick](0.4,2)--(0.8,0.4);
   			\draw[dotted][very thick](0.8,0.4)--(1.2,-0.8);
   			\draw[dotted][very thick](1.2,-0.8)--(1.6,1.6);
   			\draw[dotted][very thick](1.6,1.6)--(2,2.8);
   			\draw[dotted][very thick](2,2.8)--(2.4,-1.6);
   			\draw[dotted][very thick](2.4,-1.6)--(2.8,0);
   			
   			\draw[blue][ultra thick](0,-2)--(0,2.8);
   			\draw[blue][ultra thick](2.4,-2)--(2.4,2.8);
   			
   			\coordinate [label=above:$11$] (x) at (-1.6,-2.8);
   			\coordinate [label=above:$7$] (x) at (-1.2,-2.8);
   			\coordinate [label=above:$4$] (x) at (-0.8,-2.8);
   			\coordinate [label=above:$8$] (x) at (-0.4,-2.8);
   			\coordinate [label=above:$2$] (x) at (0,-2.8);
   			\coordinate [label=above:$10$] (x) at (0.4,-2.8);
   			\coordinate [label=above:$6$] (x) at (0.8,-2.8);
   			\coordinate [label=above:$3$] (x) at (1.2,-2.8);
   			\coordinate [label=above:$9$] (x) at (1.6,-2.8);
   			\coordinate [label=above:$12$] (x) at (2,-2.8);
   			\coordinate [label=above:$1$] (x) at (2.4,-2.8);
   			\coordinate [label=above:$5$] (x) at (2.8,-2.8);
   			
   			\coordinate [label=above:$\tau^{a}$] (x) at (-1,2.8);
   			\coordinate [label=above:$\tau^{b}$] (x) at (1.2,2.8);
   			\coordinate [label=above:$\tau^{c}$] (x) at (2.8,2.8);
   			
   			\draw[->] (3.6,0)--(5.2,0);\put(40,1.2){$\phi_{3}(\tau)$};
   			
   			\draw[step=0.4cm,gray,very thin] (5.999,-2) grid(11.2,2.8);
   			
   			\filldraw[fill=red](6.4,2.4)circle(0.1);
   			\filldraw[fill=black](6.8,0.8)circle(0.1);
   			\filldraw[fill=black](7.2,-0.4)circle(0.1);
   			\filldraw[fill=red](7.6,1.2)circle(0.1);
   			\filldraw[fill=black](8,-0.8)circle(0.1);
   			\filldraw[fill=red](8.4,2)circle(0.1);
   			\filldraw[fill=black](8.8,-1.2)circle(0.1);
   			\filldraw[fill=red](9.2,1.6)circle(0.1);
   			\filldraw[fill=red](9.6,2.8)circle(0.1);
   			\filldraw[fill=black](10,0.4)circle(0.1);
   			\filldraw[fill=black](10.4,-1.6)circle(0.1);
   			\filldraw[fill=black](10.8,0)circle(0.1);
   			
   			\draw[dotted][very thick](6.4,2.4)--(6.8,0.8);
   			\draw[dotted][very thick](6.8,0.8)--(7.2,-0.4);
   			\draw[dotted][very thick](7.2,-0.4)--(7.6,1.2);
   			\draw[dotted][very thick](7.6,1.2)--(8,-0.8);
   			\draw[dotted][very thick](8,-0.8)--(8.4,2);
   			\draw[dotted][very thick](8.4,2)--(8.8,-1.2);
   			\draw[dotted][very thick](8.8,-1.2)--(9.2,1.6);
   			\draw[dotted][very thick](9.2,1.6)--(9.6,2.8);
   			\draw[dotted][very thick](9.6,2.8)--(10,0.4);
   			\draw[dotted][very thick](10,0.4)--(10.4,-1.6);
   			\draw[dotted][very thick](10.4,-1.6)--(10.8,0);
   			
   			\draw[blue][ultra thick](8,-2)--(8,2.8);
   			\draw[blue][ultra thick](10.4,-2)--(10.4,2.8);
   			
   			\coordinate [label=above:$11$] (x) at (6.4,-2.8);
   			\coordinate [label=above:$7$] (x) at (6.8,-2.8);
   			\coordinate [label=above:$4$] (x) at (7.2,-2.8);
   			\coordinate [label=above:$8$] (x) at (7.6,-2.8);
   			\coordinate [label=above:$3$] (x) at (8,-2.8);
   			\coordinate [label=above:$10$] (x) at (8.4,-2.8);
   			\coordinate [label=above:$2$] (x) at (8.8,-2.8);
   			\coordinate [label=above:$9$] (x) at (9.2,-2.8);
   			\coordinate [label=above:$12$] (x) at (9.6,-2.8);
   			\coordinate [label=above:$6$] (x) at (10,-2.8);
   			\coordinate [label=above:$1$] (x) at (10.4,-2.8);
   			\coordinate [label=above:$5$] (x) at (10.8,-2.8);
   			
   			\coordinate [label=above:$\theta^{a}$] (x) at (7.2,2.8);
   			\coordinate [label=above:$\theta^{b}$] (x) at (9.4,2.8);
   			\coordinate [label=above:$\theta^{c}$] (x) at (11,2.8);
   		\end{tikzpicture}
   	\end{center}
   	\caption{  An example of Case 2.}\label{f4}
   \end{figure}

 From the construction of $\phi_{\ell}(\tau)$,  it is easily seen that the map $\phi_{\ell}$ preserves the relative order of the entries of $\sigma$. Hence, we have
  $\phi_{\ell}(\tau)\in S(\pi', \sigma)$, that is, the map $\phi_{\ell}$ is well-defined.

 Conversely, given any $\tau'\in S(\pi', \sigma)$, we can recover the permutation $\tau\in S(\pi, \sigma)$ as follows. If the $\ell$-th birun of $\pi'$ is increasing (resp. decreasing), then  let  $\pi'_i$ be the first (resp. last) entry of the $(\ell-2)$-th birun of $\pi'$. Suppose that $\tau'_k=\pi'_i$ for some $k\in [m+n]$.
  Then we can recover  a permutation $\tau\in S(\pi, \sigma)$ by reversing the procedure in Case 1 when the $\ell$-th birun of $\pi'$ is increasing (resp. decreasing) and $\tau'_{k+1}=\pi'_{i+1}$ (resp. $\tau'_{k-1}=\pi'_{i-1}$ ).    Otherwise, we can recover  a permutation $\tau\in S(\pi, \sigma)$ by reversing the procedure in Case 2.     So the construction
of the map  $\phi_{\ell}$ is reversible and  hence it is a bijection.

 In the following, we aim to show that ${\rm (udr, pk, des)}\tau= {\rm (udr, pk, des)}\phi_\ell(\tau)$.
 We have four  cases:  \upshape (\rmnum{1})  the $\ell$-th birun is increasing and   $\sigma^{(j-i+1)}=\emptyset$, \upshape (\rmnum{2}) the $\ell$-th birun is increasing and   $\sigma^{(j-i+1)}\neq \emptyset$, \upshape (\rmnum{3}) the $\ell$-th birun is decreasing and   $\sigma^{(j-i+1)}=\emptyset$,  and \upshape (\rmnum{4}) the $\ell$-th birun is decreasing and   $\sigma^{(j-i+1)}\neq\emptyset$.  Here we only prove the assertion for the cases  \upshape (\rmnum{1}) and \upshape (\rmnum{4}). All the other cases
can be verified by similar arguments.

 \noindent {\bf \upshape (\rmnum{1})} The $\ell$-th birun is increasing and   $\sigma^{(j-i+1)}=\emptyset$. \\
 It is easy to verify that
 $$
 {\rm des}(\tau)=
  {\rm des}(\tau^a \pi_i)+{\rm des}(\pi_{j+1}\tau^c)+ t_{\ell-1}-1+\sum_{s=1}^{j-i+1} {\rm des(\sigma^{(s)})} +\sum_{s=1}^{t_{\ell-2}-1} \delta(|\sigma^{(s)}|>0)
 $$
 and
 $$
 {\rm pk}(\tau)=
  {\rm pk}(\tau^a \pi_i)+{\rm pk}(\pi_{j+1}\tau^c)+  \sum_{s=1}^{j-i+1} {\rm epk(\sigma^{(s)})}  + \delta(|\sigma^{(t_{\ell-2} -1)}|=|\sigma^{(t_{\ell-2})}|=0).
 $$
 Here $\delta(S)=1$ if  the statement $S$ is true, and $\delta(S)=0$ otherwise.
 Similarly, we have
 $$
 {\rm des}(\phi_\ell(\tau))=
 {\rm des}(\theta^a \pi'_i)+{\rm des}(\pi'_{j+1}\theta^c)+ t'_{\ell-1}-1+\sum_{s=1}^{j-i+1} {\rm des( \sigma^{(s)})} +\sum_{s=1}^{t_{\ell-2}-1} \delta(|\sigma^{(s)}|>0),
  $$
  and
   $$
 {\rm pk}(\phi_\ell(\tau))=
  {\rm pk}(\tau^a \pi'_i)+{\rm pk}(\pi'_{j+1}\tau^c)+  \sum_{s=1}^{j-i+1} {\rm epk(\sigma^{(s)})}  + \delta(|\sigma^{(t_{\ell-2} -1)}|=|\sigma^{(t_{\ell-2})}|=0).
 $$
     As ${\rm Des}(\pi_1\pi_2\ldots \pi_{i})={\rm Des}(\pi'_1\pi'_2\ldots \pi'_{i})$ and  ${\rm Des}(\pi_{j+1}\pi_{j+2}\ldots \pi_{n})={\rm Des}(\pi'_{j+1}\pi'_{j+2}\ldots \pi'_{n})$, we have ${\rm Des}(\tau^a\pi_i)={\rm Des}(\theta^a\pi'_i)$  and ${\rm Des}(\pi_{j+1}\tau^c)={\rm Des}(\pi'_{j+1}\theta^c)$. This yields that ${\rm des}(\phi_\ell(\tau))={\rm des}( \tau)$ and ${\rm pk}(\phi_\ell(\tau))={\rm pk}( \tau)$  as $t_{\ell-1}=t'_{\ell-1}$.

    By (\ref{pkudr}),  in order to prove that ${\rm  udr}(\tau)={\rm  udr}(\phi_\ell(\tau))$, it suffices to show that $\chi^{+}(\tau)=\chi^{+}(\phi_{\ell}(\tau))$ and $\chi^{-}(\tau)=\chi^{-}(\phi_{\ell}(\tau))$.
     Assume that $\tau_x=\pi_i$ and $\tau_y=\pi_{j+1}$ for some positive integers $x$ and $y$.
     If $x=1$,  then we have $\chi^{+}(\tau)=0=\chi^{+}(\phi_{\ell}(\tau))$ since $\pi_{i}<\pi_{i+1}$ and $\pi'_{i}<\pi'_{i+1}$ guarantee that $1\notin {\rm Des}(\tau)$ and $1\notin {\rm Des}(\phi_\ell(\tau))$. If $x> 1$,   then ${\rm Des}(\tau^a\pi_i)={\rm Des}(\theta^a\pi'_i)$ implies that $\chi^{+}(\tau)=\chi^{+}(\phi_{\ell}(\tau))$.
     Notice that $\pi_{j+1}$  (resp. $\pi'_{j+1}$) is not  the last entry of the $\ell$-th birun of $\pi$ (resp. $\pi'$). This implies that $y<n+m$. Then ${\rm Des}(\pi_{j+1}\tau^c)={\rm Des}(\pi'_{j+1}\theta^c)$ implies that  $\chi^{-}(\tau)=\chi^{-}(\phi_{\ell}(\tau))$.
          So far, we have concluded that $\chi^{+}(\tau)=\chi^{+}(\phi_{\ell}(\tau))$ and $\chi^{-}(\tau)=\chi^{-}(\phi_{\ell}(\tau))$. Thus, we have ${\rm  udr}(\tau)={\rm  udr}(\phi_\ell(\tau))$ as desired.

          \noindent {\bf \upshape (\rmnum{4})} The $\ell$-th birun is deceasing and   $\sigma^{(j-i+1)}\neq \emptyset$. \\
 It is routine to check that
 $$
 {\rm des}(\tau)=
  {\rm des}(\tau^a \pi_i)+{\rm des}(\pi_{j+1}\tau^c)+ t_{\ell}-1+\sum_{s=1}^{j-i+1} {\rm des(\sigma^{(s)})} +\sum_{s=1}^{t_{\ell-1}-1} \delta(|\sigma^{(s)}|>0)
 $$
 and
 $$
 {\rm pk}(\tau)=
  {\rm pk}(\tau^a \pi_i)+{\rm pk}(\pi_{j+1}\tau^c)+  \sum_{s=1}^{j-i+1} {\rm epk(\sigma^{(s)})}  + \delta(|\sigma^{(t_{\ell-1}-1 )}|=|\sigma^{(t_{\ell-1})}|=0).
 $$
  Similarly, we have
 $$
 {\rm des}(\phi_\ell(\tau))=
 {\rm des}(\theta^a \pi'_i)+{\rm des}(\pi'_{j+1}\theta^c)+ t'_{\ell}+\sum_{s=1}^{j-i+1} {\rm des( \sigma^{(s)})} +\sum_{s=1}^{t_{\ell-1}-1} \delta(|\sigma^{(s)}|>0),
  $$
  and
   $$
 {\rm pk}(\phi_{\ell}(\tau))=
  {\rm pk}(\theta^a \pi'_i))+{\rm pk}(\pi'_{j+1}\tau^c)+  \sum_{s=1}^{j-i+1} {\rm epk(\sigma^{(s)})}  + \delta(|\sigma^{(t_{\ell-1}-1 )}|=|\sigma^{(t_{\ell-1})}|=0).
 $$
     As ${\rm Des}(\pi_1\pi_2\ldots \pi_{i})={\rm Des}(\pi'_1\pi'_2\ldots \pi'_{i})$ and  ${\rm Des}(\pi_{j+1}\pi_{j+2}\ldots \pi_{n})={\rm Des}(\pi'_{j+1}\pi'_{j+2}\ldots \pi'_{n})$, we have ${\rm Des}(\tau^a\pi_i)={\rm Des}(\theta^a\pi'_i)$  and ${\rm Des}(\pi_{j+1}\tau^c)={\rm Des}(\pi'_{j+1}\theta^c)$. This yields that ${\rm des}(\phi_\ell(\tau))={\rm des}( \tau)$ and ${\rm pk}(\phi_\ell(\tau))={\rm pk}( \tau)$  since $t'_{\ell}=t_{\ell}-1$.

    By (\ref{pkudr}),  in order to prove that ${\rm  udr}(\tau)={\rm  udr}(\phi_\ell(\tau))$, it suffices to show that $\chi^{+}(\tau)=\chi^{+}(\phi_{\ell}(\tau))$ and $\chi^{-}(\tau)=\chi^{-}(\phi_{\ell}(\tau))$.
     Assume that $\tau_x=\pi_{i}$ and $\tau_y=\pi_{j+1}$ for some positive integers $x$ and $y$. Clearly, we have $x>1$.
     Then   ${\rm Des}(\tau^a\pi_i)={\rm Des}(\theta^a\pi'_i)$ implies that $\chi^{+}(\tau)=\chi^{+}(\phi_{\ell}(\tau))$.
      If $y<n+m$, ${\rm Des}(\pi_{j+1}\tau^c)={\rm Des}(\pi'_{j+1}\theta^c)$ implies that  $\chi^{-}(\tau)=\chi^{-}(\phi_{\ell}(\tau))$. If $y=n+m$,  then we have $\chi^{-}(\tau)=0=\chi^{-}(\phi_{\ell}(\tau))$ since $\pi_{n-1}>\pi_{n}$ and $\pi'_{n-1}>\pi'_{n}$ guarantee that $n+m-1 \in {\rm Des}(\tau)$ and $n+m-1 \in {\rm Des}(\phi_\ell(\tau))$.
          So far, we have concluded that $\chi^{+}(\tau)=\chi^{+}(\phi_{\ell}(\tau))$ and $\chi^{-}(\tau)=\chi^{-}(\phi_{\ell}(\tau))$. Thus, we have ${\rm  udr}(\tau)={\rm  udr}(\phi_\ell(\tau))$ as desired.
          Hence, the map $\phi_{\ell}$ is an ${\rm(udr, pk, des)}$-preserving bijection between $S(\pi, \sigma)$ and $S(\pi', \sigma)$, completing the proof. \qed

 \noindent {\bf Proof  of Theorem \ref{th1}.} Here we only prove \upshape (\rmnum{1}).  By     similar arguments, one can verify that \upshape (\rmnum{2}), \upshape (\rmnum{3}) and \upshape (\rmnum{4}) hold. As $\pi\notin \Pi^{(1)}_{n,k,d}$, we can  find the largest integer $\ell^{(1)}$ with $
 \ell^{(1)}> 2$ such that $t_{\ell^{(1)}}\geq 3$.  Let $\pi^{(1)}$ be a permutation in $\Omega_{\ell^{(1)}}(\pi)$. By Lemma \ref{mlem1},  the map $\phi_{\ell^{(1)}}$ serves as an  ${\rm(udr, pk, des)}$-preserving bijection between $S(\pi, \sigma)$ and $S(\pi^{(1)}, \sigma)$.
 Thus we have
 $$ {(\rm udr, pk, des)}\pi= {(\rm udr, pk, des)}\pi^{(1)}    $$
 and
  $$ {(\rm udr, pk, des)}S(\pi, \sigma)= {(\rm udr, pk, des)}S(\pi^{(1)}, \sigma).    $$
  If $\pi^{(1)}\in \Pi^{(1)}_{n,k,d}$, then we stop and set $ \pi'=\pi^{(1)}$.  Otherwise,  let $t'_i$ denote the $i$-th birun of $\pi^{(1)}$.
 Then, find the largest integer $\ell^{(2)}$ with $
 \ell^{(2)}> 2$ such that $t'_{\ell^{(2)}}\geq 3$. Again  by Lemma \ref{mlem1},  the map $\phi_{\ell^{(2)}}$ serves as an  ${\rm(udr, pk, des)}$-preserving bijection between $S(\pi^{(1)}, \sigma)$ and $S(\pi^{(2)}, \sigma)$ where $\pi^{(2)}\in \Omega_{\ell^{(2)}}(\pi^{(1)})$.   We continue this process until we get some $\pi^{(s)}\in \Pi^{(1)}_{n,k,d}$. Then we set $\pi'=\pi^{(s)}$. Clearly, we have  $ {(\rm udr, pk, des)}\pi= {(\rm udr, pk, des)}\pi'    $. By   Lemma \ref{mlem1},  we have  $$ {(\rm udr, pk, des)}S(\pi, \sigma)= {(\rm udr, pk, des)}S(\pi', \sigma)  $$   as desired, completing the proof. \qed

 Now we are ready for the proof of Conjecture \ref{con1}.

 \noindent{\bf{Proof of Conjecture \ref{con1}.}}
 By  Lemma \ref{lem2},  in order to prove Conjecture \ref{con1}, it suffices to show that for any two permutations
 $\pi, \pi'\in L([n])$  with ${(\rm udr, pk, des)} \pi ={(\rm udr, pk, des)} \pi'$ and $\sigma\in L([n]+m)$, we have $ {(\rm udr, pk, des)}S(\pi, \sigma)= {(\rm udr, pk, des)}S(\pi', \sigma)    $.

 Let  $\pi, \pi'\in L([n])$  with ${(\rm udr, pk, des)} \pi ={(\rm udr, pk, des)} \pi' $ and $(\rm  pk, des)\pi=(\rm  pk, des)\pi'=(k,d)$ and let $\sigma\in L([n]+m)$.
 Notice that ${\rm Des}(\pi)={\rm Des}(\pi')$ for any permutations $\pi, \pi'\in \Pi^{(i)}_{n,k,d}$ for fixed $i\in [4]$.   Then by Lemma \ref{lem1}, we have $$ {(\rm udr, pk, des)}S(\pi, \sigma)= {(\rm udr, pk, des)}S(\pi', \sigma)   $$  when  $\pi, \pi'\in \Pi^{(i)}_{n,k,d}$ for fixed $i\in [4]$.
Otherwise,
  by Theorem  \ref{th1}, there exists two permutations $\tau, \tau'\in \Pi_{n,k,d}$ satisfying  that
  $$ {(\rm udr, pk, des)}\pi= {(\rm udr, pk, des)}\tau,   $$
  $$ {(\rm udr, pk, des)}S(\pi, \sigma)= {(\rm udr, pk, des)}S(\tau, \sigma),   $$
    $$ {(\rm udr, pk, des)}\pi'= {(\rm udr, pk, des)}\tau',   $$
   and
  $$ {(\rm udr, pk, des)}S(\pi', \sigma)= {(\rm udr, pk, des)}S(\tau', \sigma). $$
    In order to show that $ {(\rm udr, pk, des)}S(\pi, \sigma)= {(\rm udr, pk, des)}S(\pi', \sigma)    $, it remains to show that both $\tau$ and $\tau'$ are the elements of $ \Pi^{(i)}_{n,k,d}$ for some $i\in [4] $. This follows immediately from     Lemma \ref{lempk} and  the equality ${(\rm udr, pk, des)}\tau= {(\rm udr, pk, des)}\tau'$.   This completes the proof. \qed

\section{Proof of Theorem \ref{cthm}}

A  cyclic permutation statistic  ${\rm cst}$ is said to be a {\em  cyclic descent statistic} if  ${\rm cDes}([\pi])={\rm cDes}([\sigma])$ implies that ${\rm cst}([\pi])={\rm cst}([\sigma])$ for any two  cyclic permutations $[\pi]$ and $[\sigma]$.
In \cite{Doma}, Domagalski, Liang, Minnich, Sagan, Schmidt and Sietsema   derived  the following cyclic analogue of Lemma \ref{lem2}.

\begin{lemma}{\upshape   (See \cite{Doma} ,  Corollary 2.2 )}\label{lemc}
Suppose that  ${\rm cst}$ is a cyclic descent statistic. The following are equivalent.
\begin{itemize}
\item[(a)]  The statistic ${\rm cst}$ is cyclic shuffle-compatible.
\item[(b)] If  ${\rm cst} ([\pi])={\rm cst}([\pi'])$ where  $[\pi], [\pi']\in C([n])$, and $[\sigma]\in C([n]+m)$ for some $m, n\geq 1$, then $
{\rm cst}(cS([\pi], [\sigma]) )={\rm cst}(cS([\pi'], [\sigma])).
$

\end{itemize}

\end{lemma}

For any cyclic permutation $[\pi]\in C(U)$, denote by $L_i[\pi] $ the unique linear permutation in $[\pi]$ which starts with the $i$-th smallest element of $U$. For example, $L_1 [1324] =1324$, $L_2 [1324] =2413$, $L_3 [1324] =3241$ and $L_4 [1324] =4132$. It is easily seen  that for any $[\pi]\in C(U)$, we have
\begin{equation}\label{ceq1}
{\rm cdes}([\pi])={\rm cdes}(L_i[\pi] )
\end{equation}
and
\begin{equation}\label{ceq2}
{\rm cpk}([\pi])={\rm cpk}(L_i[\pi])
\end{equation}
for all $1\leq i\leq |U|$.

 For any linear permutations $\pi=\pi_1\pi_2\ldots \pi_n\in L([n])$ and $\sigma\in L([n]+m)$,  denote by $S'(\pi,\sigma)$ the set of permutations $\tau=\tau_1\tau_2\ldots \tau_{n+m}\in S(\pi, \sigma)$ with $\tau_1=\pi_1$ and $\tau_{n+m}=\pi_n$.  Denote by $L'(U)$ the set of linear permutations $\pi\in L(U)$ which start with the   smallest element of $U$ and end  with the second smallest element of $U$.

The following theorem  will play an essential role in the proof of Theorem  \ref{cthm}.
\begin{theorem}\label{cth1}
Let $\pi$ and $ \pi'$ be  permutations in $ L'([n])$ with $({\rm pk, des }) \pi =({\rm pk, des  }) \pi' $ and let $\sigma\in L([n]+m)$ for some $m \geq 1$,  $n>2$,  and $k,d\geq 0$.   Then we have
$$
(\rm pk, des) S'(\pi, \sigma) =(\rm pk, des) S'(\pi', \sigma) .
$$
\end{theorem}

Before we prove Theorem \ref{cth1}, we need  the following two lemmas.

 \begin{lemma}\label{cmlem1}
Let $\pi\in L'([n])$ be a permutation  with  ${\rm type}(\pi)=(t_1, t_2, \ldots, t_{2k})$ such that $t_{\ell}\geq 3$ for some $\ell\geq 3$ and let $\sigma\in L([n]+m)$ for some $m \geq 1$ and $n>2$. The map $\phi_\ell$ induces a   ${\rm ( pk, des)}$-preserving  bijection between   $S'(\pi, \sigma)$  and $ S'(\pi', \sigma)$ for any permutation $\pi'\in \Omega_{\ell}(\pi)\cap L'([n])$.

\end{lemma}
\pf  From the construction of the map $\phi_{\ell}$, one can easily check that for any $\tau\in S'(\pi, \sigma)$, we have $\phi_{\ell}(\tau)\in S'(\pi', \sigma) $ as desired, completing the proof.  \qed

\begin{lemma}\label{cmlem2}
Let $\pi\in L'([n])$ be a permutation  with $({\rm pk, des }) \pi =(k,d)$ and let $\sigma\in L([n]+m)$ for some $m \geq 1$, $n>2$,  and $k,d\geq 0$. If $\pi\notin  \Pi^{(1)}_{n,k,d}$, then there exists  a permutation $\pi'\in  \Pi^{(1)}_{n,k,d}\cap L'([n])$   such that  $$ {(\rm   pk, des)}S'(\pi, \sigma)= {(\rm  pk, des)}S'(\pi', \sigma).  $$
 \end{lemma}
 \pf
 Since $\pi\in L'([n])$, we have ${\rm type}(\pi)=(t_1, t_2, \ldots, t_{2k})$  and $\chi^{+}(\pi)=0$.
As $\pi\notin \Pi^{(1)}_{n,k,d}$, we can  find the largest integer $\ell^{(1)}$ with $
 \ell^{(1)}> 2$ such that $t_{\ell^{(1)}}\geq 3$.  Let $\pi^{(1)}$ be a permutation in $\Omega_{\ell^{(1)}}(\pi)\cap L'([n])$. By Lemma \ref{cmlem1},  the map $\phi_{\ell^{(1)}}$ serves as a   ${\rm(  pk, des )}$-preserving bijection between $S'(\pi, \sigma)$ and $S'(\pi^{(1)}, \sigma)$.
  If $\pi^{(1)}\in \Pi^{(1)}_{n,k,d}$, then we set $\pi'=\pi^{(1)}$.  Otherwise,  let $t'_i$ denote the $i$-th birun of $\pi^{(1)}$.
 Then, find the largest integer $\ell^{(2)}$ with $
 \ell^{(2)}> 2$ such that $t'_{\ell^{(2)}}\geq 3$. Again  by Lemma \ref{cmlem1},  the map $\phi_{\ell^{(2)}}$ serves as a   ${\rm(pk, des)}$-preserving bijection between $S'(\pi^{(1)}, \sigma)$ and $S'(\pi^{(2)}, \sigma)$ where $\pi^{(2)}\in \Omega_{\ell^{(2)}}(\pi^{(1)})\cap L'([n])$.      We continue this process until we get some   $\pi^{(s)}\in \Pi^{(1)}_{n,k,d}\cap L'([n])$. Let $\pi'=\pi^{(s)}$.   By   Lemma \ref{cmlem1},  we have  $$ {(\rm  pk, des)}S'(\pi, \sigma)= {(\rm  pk, des)}S'(\pi', \sigma)  $$   as desired, completing the proof.  \qed

\noindent{\bf Proof of Theorem \ref{cth1}.} Assume that$({\rm pk, des }) \pi =({\rm pk, des  }) \pi'=(k,d) $.  If $\pi, \pi'\in \Pi^{(1)}_{n,k,d}\cap L'([n])$, we first  describe  a map $\psi: S' (\pi, \sigma)\rightarrow S' (\pi', \sigma) $ as follows.   For any $\tau\in S'(\pi, \sigma)$, define $ \psi(\tau)$ to be the permutation obtained from $\tau$ by replacing each $\pi_i$ by $\pi'_i$ for all $1\leq i\leq n$. Clearly, we have  $\psi(\tau)\in S'(\pi', \sigma)$ and ${\rm Des}(\tau)={\rm Des}(\psi(\tau))$, which implies that  $ {(\rm   pk, des)}(\tau)={(\rm   pk, des)}(\psi(\tau))$. Clearly, the map $\psi$ is reversible  and hence it is a bijection.
Therefore, we have
$$
(\rm pk, des) S'(\pi, \sigma) =(\rm pk, des) S'(\pi', \sigma) .
$$ when $\pi, \pi'\in \Pi^{(1)}_{n,k,d}\cap L'([n])$.
Otherwise,
  by Lemma  \ref{cmlem2}, there exists two permutations $\tau, \tau'\in \Pi^{(1)}_{n,k,d}\cap L'([n])$ satisfying  that
  $$ {(\rm   pk, des)}S'(\pi, \sigma)= {(\rm   pk, des)}S'(\tau, \sigma),   $$
    and
  $$ {(\rm   pk, des)}S'(\pi', \sigma)= {(\rm  pk, des)}S'(\tau', \sigma). $$
   Then  the equality $ {(\rm  pk, des)}S'(\pi, \sigma)= {(\rm   pk, des)}S'(\pi', \sigma)    $  follows  immediately form the equality $$
(\rm pk, des) S'(\tau, \sigma) =(\rm pk, des) S'(\tau', \sigma) .
$$     This completes the proof. \qed

Now we are ready for the proof of Theorem \ref{cthm}.

 \noindent{\bf{Proof of Theorem  \ref{cthm}.}}
 For any $[\pi]\in C[n]$ for $n\geq 2$, let $L^{+}_1[\pi]$ denote the permutation obtained from $L_1[\pi] $ by increasing  each element of $[n]\setminus \{1\}$ by one and inserting a $2$ at the end of $L_1[\pi]$. For example, $L^{+}_1[1342]=14532$.  Clearly, we have $L^{+}_1[\pi]\in L'([n+1])$.  It is easily seen that $({\rm  cpk,  cdes} ) L_1[\pi]=({\rm  pk,  des} ) L^{+}_1[\pi]$.
 For any two cyclic  permutations
 $[\pi], [\pi']\in C([n]) $  with ${(\rm   cpk, cdes)} [\pi] ={(\rm   cpk, cdes)} [\pi']$ and $\sigma\in C([n]+m)$.
 Then, we have
 $$
\begin{array}{lll}
{ \rm   (cpk, cdes) } cS([\pi], [\sigma])&=& \bigcup\limits_{[\tau]\in cS([\pi], [\sigma])}\{ { \rm   (cpk, cdes) }[\tau]\}\\
&&\\
&=& \bigcup\limits_{[\tau]\in cS([\pi], [\sigma])}\{ { \rm   (cpk, cdes) }L_1[\tau]\}\,\,\,\,\,\,\, \mbox{(by (\ref{ceq1}) and (\ref{ceq2})   )}\\
&&\\
&=&    \bigcup\limits_{[\tau]\in  cS( [\pi] , [\sigma])}\{ { \rm   (pk, des) }L^{+}_1[\tau]\}.
\end{array}
$$
It is easy to check that
$$
\{L^{+}_1[\tau] \mid  [\tau]\in  cS( [\pi] , [\sigma])\}=\bigcup\limits_{j=1}^m S'(L^+_1[\pi] , L_j[\sigma]).
$$
Hence, we have
\begin{equation}\label{ceq3}
{ \rm   (cpk, cdes) } cS([\pi], [\sigma])=\bigcup\limits_{j=1}^{m} \bigcup\limits_{\tau\in S'(L^+_1[\pi] , L_j[\sigma])} \{{ \rm   (pk, des) }\tau\}.
 \end{equation}

 Let $m\geq 1$ and $n\geq 2$.
 By  Lemma \ref{lemc},  in order to prove Theorem  \ref{cthm}, it suffices to show that for any two cyclic  permutations
 $[\pi], [\pi']\in C([n])$   with ${(\rm   cpk, cdes)} [\pi] ={(\rm   cpk, cdes)} [\pi']$ and $\sigma\in C([n]+m)$, we have $$ {(\rm   cpk, cdes)}cS([\pi], [\sigma])= {(\rm   cpk, cdes)}cS([\pi'], [\sigma]).   $$
As ${(\rm   cpk, cdes)} [\pi] ={(\rm   cpk, cdes)} [\pi']$, we have ${(\rm   pk, des)} L^+_1[\pi] ={(\rm   pk, des)} L^+_1[\pi']$. Then by Theorem \ref{cth1}, we deduce that
\begin{equation}\label{ceq4}
   \bigcup\limits_{\tau\in S'(L^+_1[\pi] , L_j[\sigma])} \{{ \rm   (pk, des) }\tau\}=\bigcup\limits_{\tau\in S'(L^+_1[\pi'] , L_j[\sigma])} \{{ \rm   (pk, des) }\tau\}
 \end{equation}
 for all $1\leq j\leq m$.
Combining (\ref{ceq3}) and (\ref{ceq4}), we have $${ \rm   (cpk, cdes) } cS([\pi], [\sigma])={ \rm   (cpk, cdes) } cS([\pi'], [\sigma])$$ as desired, completing the proof.
  \qed

\vspace{0.5cm}
 \noindent{\bf Acknowledgments.}
 This work was supported by  the National Natural Science Foundation of China (12071440).


\begin{thebibliography}{100}

\bibitem{Adin}
R.  M. Adin, I.  M. Gessel, V.  Reiner, Y.  Roichman, Cyclic quasi-symmetric functions,
{\em Israel J. Math.},  {\bf 243} (2021),  437--500.

\bibitem{Baker}
D. Baker-Jarvis, B.E. Sagan, Bijective proofs of shuffle compatibility results, {\em Adv.  Appl. Math.}, {\bf 113} (2020),  101973.

\bibitem{Bona}
 M. B\'{o}na, Combinatorics of Permutations. CRC Press, 2004




 \bibitem{Doma}
 R. Domagalski, J. Liang, Q. Minnich,  B.E. Sagan, J. Schmidt, A. Sietsema, Cyclic shuffle compatibility,
{\em S\'em. Lothar. Combin.},  {\bf 85} ([2020--2021]), Art. B85d, 11 pp.

\bibitem{Gar}
A. M. Garsia and I. M. Gessel, Permutation statistics and partitions, {\em Adv. in Math.},  {\bf 31} (1979),  288--305.

\bibitem{Gessel}
I.M. Gessel, Y. Zhuang, shuffle-comapatible permutation statitics, {\em Adv. Math.}, {\bf 332} (2018), 85--141.

\bibitem{Goulden}
I. P. Goulden, A bijective proof of Stanley's shuffling theorem, {\em Trans. Amer. Math. Soc.}, { \bf 288}
(1985),  147--160.


\bibitem{Grinberg}
D. Grinberg, Shuffle-compatible permutation statistics II: the exterior peak set, {\em Electron. J.
Combin.},  {\bf 25} (2018),  P4.17.

\bibitem{Ji}
K.Q. Ji, D.T.X. Zhang, Stanley's shuffle theorem and insertion lemma, arXiv:2203.13543v1 [math.CO].

\bibitem{Ji2}
K.Q. Ji, D.T.X. Zhang, A cyclic analogue of Stanley's shuffle theorem, arXiv:2205.03188v1 [math.CO].

\bibitem{Novick}
M. Novick, A bijective proof of a major index theorem of Garsia and Gessel, {\em Electron. J.
Combin.},  {\bf 17 (1)} (2010),  64.


\bibitem{Oguz}
 E. K.  O$\breve{g}$uz, A counter example to the shuffle compatibility conjecture, arXiv:1807.01398
[math.CO].

\bibitem{Stadler}
J. D. Stadler, Stanley's shuffling theorem revisited, {\em J. Combin. Theory Ser. A },   {\bf 88}  (1999),
176--187.


\bibitem{Stanley}
R.P. Stanley, Ordered Structures and Partitions, Memoirs of the American Mathematical
Society, vol. 119, American Mathematical Society, Providence, R.I., 1972
	
\end{thebibliography}
\end{document}